\documentclass[12pt]{article}
\usepackage{amssymb}
\usepackage{amsfonts}
\usepackage{latexsym}
\usepackage{amsthm}

\begin{document}


\renewcommand{\theequation}{\arabic{section}.\arabic{equation}}
\def\[{\begin{equation}}              \def\]{\end{equation}}
\def\lb{\label}                       \newcommand{\er}[1]{{\rm(\ref{#1})}}
\theoremstyle{plain}
\newtheorem{theorem}{\bf Theorem}[section]
\newtheorem{lemma}[theorem]{\bf Lemma}
\newtheorem{corollary}[theorem]{\bf Corollary}
\newtheorem{proposition}[theorem]{\bf Proposition}
\newtheorem{definition}[theorem]{\bf Definition}
\newtheorem{remark}[theorem]{\it Remark}

\def\a{\alpha}  \def\cA{{\cal A}}     \def\bA{{\bf A}}  \def\mA{{\mathscr A}}
\def\b{\beta}   \def\cB{{\cal B}}     \def\bB{{\bf B}}  \def\mB{{\mathscr B}}
\def\g{\gamma}  \def\cC{{\cal C}}     \def\bC{{\bf C}}  \def\mC{{\mathscr C}}
\def\G{\Gamma}  \def\cD{{\cal D}}     \def\bD{{\bf D}}  \def\mD{{\mathscr D}}
\def\d{\delta}  \def\cE{{\cal E}}     \def\bE{{\bf E}}  \def\mE{{\mathscr E}}
\def\D{\Delta}  \def\cF{{\cal F}}     \def\bF{{\bf F}}  \def\mF{{\mathscr F}}
\def\c{\chi}    \def\cG{{\cal G}}     \def\bG{{\bf G}}  \def\mG{{\mathscr G}}
\def\z{\zeta}   \def\cH{{\cal H}}     \def\bH{{\bf H}}  \def\mH{{\mathscr H}}
\def\e{\eta}    \def\cI{{\cal I}}     \def\bI{{\bf I}}  \def\mI{{\mathscr I}}
\def\p{\psi}    \def\cJ{{\cal J}}     \def\bJ{{\bf J}}  \def\mJ{{\mathscr J}}
\def\vT{\Theta} \def\cK{{\cal K}}     \def\bK{{\bf K}}  \def\mK{{\mathscr K}}
\def\k{\kappa}  \def\cL{{\cal L}}     \def\bL{{\bf L}}  \def\mL{{\mathscr L}}
\def\l{\lambda} \def\cM{{\cal M}}     \def\bM{{\bf M}}  \def\mM{{\mathscr M}}
\def\L{\Lambda} \def\cN{{\cal N}}     \def\bN{{\bf N}}  \def\mN{{\mathscr N}}
\def\m{\mu}     \def\cO{{\cal O}}     \def\bO{{\bf O}}  \def\mO{{\mathscr O}}
\def\n{\nu}     \def\cP{{\cal P}}     \def\bP{{\bf P}}  \def\mP{{\mathscr P}}
\def\r{\rho}    \def\cQ{{\cal Q}}     \def\bQ{{\bf Q}}  \def\mQ{{\mathscr Q}}
\def\s{\sigma}  \def\cR{{\cal R}}     \def\bR{{\bf R}}  \def\mR{{\mathscr R}}
\def\S{\Sigma}  \def\cS{{\cal S}}     \def\bS{{\bf S}}  \def\mS{{\mathscr S}}
\def\t{\tau}    \def\cT{{\cal T}}     \def\bT{{\bf T}}  \def\mT{{\mathscr T}}
\def\f{\phi}    \def\cU{{\cal U}}     \def\bU{{\bf U}}  \def\mU{{\mathscr U}}
\def\F{\Phi}    \def\cV{{\cal V}}     \def\bV{{\bf V}}  \def\mV{{\mathscr V}}
\def\P{\Psi}    \def\cW{{\cal W}}     \def\bW{{\bf W}}  \def\mW{{\mathscr W}}
\def\o{\omega}  \def\cX{{\cal X}}     \def\bX{{\bf X}}  \def\mX{{\mathscr X}}
\def\x{\xi}     \def\cY{{\cal Y}}     \def\bY{{\bf Y}}  \def\mY{{\mathscr Y}}
\def\X{\Xi}     \def\cZ{{\cal Z}}     \def\bZ{{\bf Z}}  \def\mZ{{\mathscr Z}}
\def\O{\Omega}

\def\ve{\varepsilon}   \def\vt{\vartheta}    \def\vp{\varphi}    \def\vk{\varkappa}

\def\Z{{\mathbb Z}}    \def\R{{\mathbb R}}   \def\C{{\mathbb C}}
\def\T{{\mathbb T}}    \def\N{{\mathbb N}}   \def\dD{{\mathbb D}}


\def\la{\leftarrow}              \def\ra{\rightarrow}            \def\Ra{\Rightarrow}
\def\ua{\uparrow}                \def\da{\downarrow}
\def\lra{\leftrightarrow}        \def\Lra{\Leftrightarrow}


\let\ge\geqslant                 \let\le\leqslant
\def\/{\over}                    \def\iy{\infty}
\def\sm{\setminus}               \def\es{\emptyset}
\def\ss{\subset}                 \def\ts{\times}
\def\pa{\partial}                \def\os{\oplus}
\def\ev{\equiv}                  \def\1{1\!\!1}
\def\iint{\int\!\!\!\int}        \def\iintt{\mathop{\int\!\!\int\!\!\dots\!\!\int}\limits}
\def\el2{\ell^{\,2}}


\def\qq{\quad}                   \def\qqq{\qquad}
\def\lt{\biggl}                  \def\rt{\biggr}
\def\no{\noindent}               \def\ol{\overline}
\def\wt{\widetilde}              \def\wh{\widehat}
\newcommand{\nt}[1]{{\mathop{#1}\limits^{{}_{\,\bf{\sim}}}}\vphantom{#1}}
\newcommand{\nh}[1]{{\mathop{#1}\limits^{{}_{\,\bf{\wedge}}}}\vphantom{#1}}
\newcommand{\nc}[1]{{\mathop{#1}\limits^{{}_{\,\bf{\vee}}}}\vphantom{#1}}
\newcommand{\oo}[1]{{\mathop{#1}\limits^{\,\circ}}\vphantom{#1}}
\newcommand{\po}[1]{{\mathop{#1}\limits^{\phantom{\circ}}}\vphantom{#1}}


\def\Im{\mathop{\rm Im}\nolimits}
\def\Iso{\mathop{\rm Iso}\nolimits}
\def\Ker{\mathop{\rm Ker}\nolimits}
\def\Ran{\mathop{\rm Ran}\nolimits}
\def\Re{\mathop{\rm Re}\nolimits}
\def\Tr{\mathop{\rm Tr}\nolimits}
\def\arg{\mathop{\rm arg}\nolimits}
\def\const{\mathop{\rm const}\nolimits}
\def\det{\mathop{\rm det}\nolimits}
\def\diag{\mathop{\rm diag}\nolimits}
\def\dim{\mathop{\rm dim}\nolimits}
\def\dist{\mathop{\rm dist}\nolimits}
\def\rank{\mathop{\rm rank}\limits}
\def\res{\mathop{\rm res}\limits}
\def\sign{\mathop{\rm sign}\nolimits}
\def\span{\mathop{\rm span}\nolimits}
\def\supp{\mathop{\rm supp}\nolimits}

\newcommand\matr[4]{\left(\begin{array}{cc}#1 & #2 \cr #3 & #4\end{array}\right)}
\newcommand\vect[2]{\left(\begin{array}{c} #1 \cr #2 \end{array}\right)}

\title {Parametrization of the isospectral set
for the vector-valued Sturm-Liouville problem}

\author{Dmitry Chelkak\begin{footnote}
{Dept. of Math. Analysis, Math. Mech. Faculty, St.Petersburg State University. Universitetskij
pr. 28, Staryj Petergof, 198504 St.Petersburg, Russia, e-mail: delta4@math.spbu.ru }
\end{footnote}
and Evgeny Korotyaev\begin{footnote} { Institut f\"ur  Mathematik,  Humboldt Universit\"at zu
Berlin, Rudower Chaussee 25, 12489, Berlin, Germany, e-mail: evgeny@math.hu-berlin.de}
\end{footnote}
}

\maketitle

\begin{abstract}
\no We obtain a parametrization of the isospectral set of  matrix-valued potentials for the
vector-valued Sturm-Liouville problem on a finite interval.

\vskip 2pt

\no {\bf AMS Classification:} 34A55, 34B24.

\vskip 2pt

\no {\bf Keywords:} Sturm-Liouville problem, matrix potentials, isospectral potentials.
\end{abstract}

\section{Introduction and main results}

Consider the inverse problems for the self-adjoint operator $H$ in $L^2(0,1)^N$, given  by
\[
\label{Hdef} H\p=-\p''+V(x)\p=\l\p,\qquad \p(0)=\p(1)=0,
\]
where $V=\!V^*\!\in\!L^1(0,1)$ is some $N\ts N$ matrix-valued potential and $\p$ is a
vector-valued function. Let
$$
\l_1<\l_2<\dots<\l_\a<\dots
$$
be the eigenvalues of $H$, where each $\l_\a$, $\a\ge 1$,  has multiplicity $k_\a\in [1,N]$,
i.e. $k_\a$ is the number of eigenfunctions corresponding to the eigenvalue $\l_\a$.

In the scalar case this problem is well-known, including the complete characterization of the
set of spectral data (eigenvalues and norming constants) that correspond to various classes of
potentials (see \cite{L}, \cite{M}, \cite{PT}). On the contrary, in the matrix case only some
particular results are known. Thus, it is known that the matrix-valued Weyl-Titchmarsh
function (see (\ref{mDef})) uniquely determines the potential $V$ (see \cite{Mal} or
\cite{Yu}, where this Borg-type uniqueness theorem was proved in two different ways). Jodeit
and Levitan \cite{JL1} and Chern \cite{Ch} have constructed some isospectral sets of
potentials, which not only have the same spectrum as $V$ but also the same initial data
(boundary values) of eigenfunctions. Some other results about the inverse spectral problem for
the vector-valued Schr\"odinger equation were obtained in \cite{CHGL}, \cite{CK}, \cite{Ca},
\cite{ChSh}, \cite{JL2}, \cite{Sh}, \cite{SP}.

We denote by $\vp(x,\l)$ the $N\ts N$ matrix-valued solution of the equation
$-\vp''+V\vp=\l\vp$ such that $\vp(0,\l)=0$ and $\vp'(0,\l)=I$, where $I=I_N$ is the identity
matrix. Introduce the matrices
$$
S_\a=\int_0^1(\vp^*\vp)(t,\l_\a)dt=S_\a^*>0,\quad \a\ge 1.
$$

\begin{definition}[{\bf Spectral data}]
For each eigenvalue $\l_\a$, $\a\ge 1$, we define {\bf the subspace}
\[
\label{cEDef} {\cE_\a}=\Ker \vp(1,\l_\a)=
\left\{h\in\C^N:\vp(1,\l_\a)h=0\right\}\ss\C^N,\qquad \dim\cE_\a=k_\a,
\]
the orthogonal projector $P_\a:\C^N\!\to\! \cE_\a$ and {\bf the positive self-adjoint
operator} $g_\a:\cE_\a\!\to\!\cE_\a$ given by
$$
{g_\a}= G_\a\Big|_{\cE_\a},\quad {where}\quad G_\a=P_\a S_\a P_\a.
$$
\end{definition}
\no {\it Remark.} Each solution of the equation $-\p''+V\p=\l\p$ under the condition $\p(0)=0$
has the form $\vp(x,\l)h$ for some $h\in \C^N$. Then, each eigenfunction $\p_\a$ such that
$H\p_\a=\l\p_\a$  has the form $\p_\a(x)=\vp(x,\l_\a)h$ for some $h\in \cE_\a$\,. In
particular, we get $\dim \cE_\a=k_\a$. Moreover, the following identity is fulfilled:
$$
\int_0^1\|\vp(x,\l_\a)h\|^2dx= \langle h, g_\a h\rangle,\quad h\in \cE_\a\,,
$$
where $\langle u,v\rangle= u^*v$ is the scalar product of two vectors and $\|u\|^2=\langle
u,u\rangle$.

Below we will sometimes write $\vp(x,\l,V), \l_\a(V),\dots$, instead of $\vp(x,\l),
\l_\a,\dots$, when several potentials are being dealt with.

Let $\wt{\l}_\a=\l_\a(\wt{V})$, $\wt{\vp}(x,\l)=\vp(x,\l,\wt{V})$ and so on. For $V=V^*\in
L^1(0,1)$ we introduce the isospectral set of potentials by
\[
\label{IsoDef}
\Iso(V)=\left\{\wt V=\wt V^*\in L^1(0,1):\wt{\l}_\a=\l_\a,\ \wt{k}_\a=k_\a \ {\rm
for\ all }\ \a\ge 1\right\}.
\]

Our goal is to show that the spectral data $\{(\cE_\a,g_\a)\}_{\a\ge 1}$ give the "proper"\
parametrization of the set $\Iso(V)$. For the sake of the reader, we start with the simple
result that the parameters $\{(\l_\a,\cE_\a,g_\a)\}_{\a\ge 1}$  determine the potential
uniquely (see Theorem \ref{*UniqThm}). In fact, these data are closely related to the residues
of the Weyl-Titchmarsh function (see Proposition \ref{*ResM=}). Our main result is Theorem
\ref{*EGTransThm} which shows that $\{(\cE_\a,g_\a)\}_{\a\ge 1}$ are {\bf free parameters}.
Namely, we prove that each $\cE_\a$ and $g_\a$ can be changed in an almost arbitrary way, when
all other parameters $\{(\cE_\b,g_\b)\}_{\b\ne\a}$ and the spectrum are fixed.

\begin{theorem}[{\bf Uniqueness}]
\label{*UniqThm} Let $\wt{V}\in \Iso(V)$ for some $V=V^*\in L^1(0,1)$  and let
$\wt{\cE}_\a=\cE_\a$ and $\wt{g}_\a=g_\a$ for all $\a\ge 1$. Then $\wt{V}=V$.
\end{theorem}

Introduce  the {\bf "forbidden subspace"}
\[
\label{FaDef} \cF_\a=\C^N\ominus (S_\a(\cE_\a)),\quad \a\ge 1,
\]
where $\C^N\ominus \cS = \{v\in\C^N:\langle v,u \rangle=0\ {\rm for\ all}\ u\in\cS\}$ is the
orthogonal subspace to $\cS\ss\C^N$. Note that $\dim\cF_\a=N\!-\!k_\a$, since
$\dim\cE_\a=k_\a$ and $S_\a>0$. We formulate our main result.

\begin{theorem}
\label{*EGTransThm} Let $V=V^*\in L^1(0,1)$ and $\a\ge 1$. Then the mapping
$\Phi_\a:\wt{V}\mapsto (\wt{\cE}_\a,\wt{g}_\a)$ is a {bijection} between the set of potentials
$$
\left\{\wt{V}\in\Iso(V):\ \wt{\cE}_\b=\cE_\b,\ \wt{g}_\b=g_\b\ { for\ all}\ \b\ne\a\right\}
$$
and the following set of pairs $(\wt{\cE}_\a,\wt{g}_\a)$:
$$
\left\{({\cE},{g}):\ {g}={g}^*>0\ { is\ an\ operator\ in\ some\ subspace\ }{\cE}\ss\C^N\right.
$$
$$
\left. \ {with\ }\dim\cE=k_\a\ {such\ that\ } {\cE}\cap\cF_\a=\{0\}\right\}.
$$
\end{theorem}
\no {\it Remark.\ } i) Substituting $\wt{V}=V$, we obtain $\cE_\a\cap\cF_\a=\{0\}$.

\no ii) If $\wt{\cE}_\a=\cE_\a$, then there are no restrictions on the changing of the
"norming matrix"\ $g_\a$. This case is similar to the scalar case and \cite{JL1}, \cite{Ch}.

\no iii) If we change $\wt{\cE}_\a$, then there exists only one restriction
$\wt{\cE}_\a\cap\cF_\a=\{0\}$. Such an effect is absent in the scalar case. As far as the
authors know, this is the first result in this direction.

\no iv) Note that we use an explicit procedure (see Theorem \ref{*TransfThm}), which is based
on the so-called Darboux transform. Therefore, the result of any finite number of such changes
can be expressed explicitly in terms of the initial potential.

The next Proposition shows that each "forbidden subspace"\ $\cF_\a$, $\a\ge 1$, doesn't depend
on the "norming matrices" $\{g_\b\}_{\b\ge 1}$. Namely, it is uniquely determined by the
spectrum and all subspaces $\{\cE_\b\}_{\b\ne\a}$.

\begin{proposition}
\label{*FaProp} Let $\wt{V}\in\Iso(V)$ for some $V=V^*\in L^1(0,1)$. Fix some $\a\ge 1$ and
let $\wt{\cE}_\b=\cE_\b$ for all $\b\ne\a$. Then $\wt{\cF}_\a=\cF_\a$ and
$\wt{\cE}_\a\cap\cF_\a=\{0\}$.
\end{proposition}

In order to illustrate the "forbidden subspaces", we give the following simple example.

\begin{proposition}[{\bf Example}]
\label{*ExProp} Let $N=2$ and $V=V^*\in L^1(0,1)$ be such that $k_1=k_2=1$ and $k_\a=2$ for
all $\a\ge 3$. Then $\cE_1\cap\cE_2=\{0\}$, $\cF_1=\cE_2$ and $\cF_2=\cE_1$.
\end{proposition}

Finally, we give the connection between our spectral data and the matrix-valued
Weyl-Titchmarsh function $m(\l)$ given by
\[
\label{mDef} m(\l)=(\c'\c^{-1})(0,\l),\quad \l\in\C,
\]
where $\c(x,\l)$ is the matrix-valued solution of the equation $-\c''+V\c=\l\c$ such that
$\c(1,\l)=0$ and $\c'(1,\l)=I$.

\begin{proposition}
\label{*ResM=} Let $V=V^*\in L^1(0,1)$. Then the function $m$ is analytic in $\C\sm
\bigcup_{\a\ge 1}\{\l_\a\}$ and satisfies the identity $m(\l)=m^*(\ol{\l})$. Moreover, each
point $\l_\a$, $\a\ge 1$, is a simple pole of $m$  and
$$
\res_{\l=\l_\a}m(\l)\Big|_{\cE_\a}\!=\,-g_\a^{-1},\qquad
\res_{\l=\l_\a}m(\l)\Big|_{\C^N\ominus\cE_\a}\!\!=\,\,0,\quad \a\ge 1.
$$
\end{proposition}

We describe the plan of the paper. In Sect. 2 we prove some preliminary Lemmas, Theorem
\ref{*UniqThm} and Proposition \ref{*ResM=}. In Sect. 3 we prove Theorem \ref{*EGTransThm} and
Propositions \ref{*FaProp},\ref{*ExProp}.

\section{Preliminaries}
\setcounter{equation}{0}

Repeating the standard arguments (see \cite{PT}, p. 13--15), we obtain the following
asymptotics:
\[
\label{VpAsympt} \vp(x,z^2)=\frac{\sin zx}{z}\,\cdot I -\frac{\cos
zx}{2z^2}\int_0^xV(t)dt+o\lt(\frac{e^{|\Im z|x}}{|z|^2}\rt),
\]
\[
\label{Vp'Asympt} \vp'(x,z^2)= \cos zx \cdot I+\frac{\sin zx}{2z}\int_0^xV(t)dt+
o\lt(\frac{e^{|\Im z|x}}{|z|}\rt),
\]
as $|z|\to\infty$ for all $x\in [0,1]$ and $V\in L^1(0,1)$. Also, note that
\[
\label{ChiIdentity} \c(x,\l,V)=-\vp(1\!-\!x,\l,V^\sharp),\quad {\rm  where}\quad V^\sharp
(t)=V(1-t),\ t\in [0,1].
\]

\begin{lemma}
\label{*Phi2Ident} Let $V=V^*\in L^1(0,1)$. Then

\no (i) for each $(x,\l)\in [0,1]\ts\C$ the following identities are fulfilled:
\[
\label{VpIdent} \vp^*(x,\ol{\l})\vp'(x,\l)=(\vp')^*(x,\ol{\l})\vp(x,\l),\qquad
\c^*(x,\ol{\l})\c'(x,\l)=(\c')^*(x,\ol{\l})\c(x,\l),
\]
\[
\label{VpCIdent}
\c^*(0,\ol{\l})=\c^*(x,\ol{\l})\vp'(x,\l)-(\c')^*(x,\ol{\l})\vp(x,\l)=-\vp(1,\l).
\]
\no (ii) for each $\a\ge 1$ the following identity is fulfilled:
\[
\label{LIdent} G_\a=P_\a\,[\dot{\vp}^*\vp'](1,\l_\a)P_\a.
\]
\end{lemma}
\begin{proof}
(i) The function $\e(x)=\c^*(x,\ol{\l})\vp'(x,\l)-(\c')^*(x,\ol{\l})\vp(x,\l)$ satisfies the
equation
$$
\e'(x)=\c^*(x,\ol{\l})\left((V(x)\!-\!\l I)-(V^*(x)\!-\!\ol\l I)^*\right)\vp(x,\l)=0,\qquad
x\in [0,1].
$$
Due to $\e(0)=\c^*(0,\ol{\l})$ and $\e(1)=-\vp(1,\l)$, we obtain (\ref{VpCIdent}). The proof
of (\ref{VpIdent}) is similar.

\no (ii) Note that $-\dot\vp''=(\l I-V)\dot\vp+\vp$. This gives
$\left[\dot{\vp}^*\vp'-(\dot{\vp}')^*\vp\right]'(x,\l_\a)=[\vp^*\vp](x,\l_\a)$, since
$\l_\a\in\R$. Therefore,
$$
G_\a= P_\a \left[\dot{\vp}^*\vp'-(\dot{\vp}')^*\vp\right](1,\l_\a) P_\a=
P_\a[\dot{\vp}^*\vp'](1,\l_\a)P_\a,
$$
where we have used $\vp(1,\l_\a)P_\a=0$.
\end{proof}

Introduce the matrices
\[
\label{DaDef} Z_\a=\dot\vp(1,\l_\a)P_\a +\vp(1,\l_\a)P_\a^\bot,\quad {\rm where} \quad
P_\a^\bot=I-P_\a, \quad \a\ge 1.
\]
\begin{lemma}
\label{*ZaProp} Let $V\!=\!V^*\!\in\!L^1(0,1)$. Then (i) $\det Z_\a\ne 0$ for all $\a\ge 1$.

\no (ii) Each $\l_\a$, $\a\!\ge\!1$, is the root of the entire function $\det\vp(1,\l)$ of the
multiplicity \nolinebreak $k_\a$. The function $\det\vp(1,\l)$ has no other roots. Moreover,
the following asymptotics is fulfilled:
\[
\label{Vp-1=}
\vp^{-1}(1,\l)=((\l\!-\!\l_\a)^{-1}P_\a+P_\a^\bot)(Z_\a^{-1}+O(\l\!-\!\l_\a))\quad {\rm as}\ \
\l\to\l_\a.
\]

\no (iii) Let $\xi(\l)$, $\l\in\C$, be the entire $N\!\times\!N$ matrix-valued function such
that $\xi(\l_\a)P_\a=0$ for all $\a\ge 1$. Then, $\xi(\l)\vp^{-1}(1,\l)$ is the entire
matrix-valued function.
\end{lemma}

\begin{proof}
(i) Suppose that $\dot\vp(1,\l_\a)P_\a h+\vp(1,\l_\a)P_\a^\bot h=0$ for some vector $h\in
\C^N$. Using (\ref{LIdent}) and (\ref{VpIdent}), we obtain
$$
\langle P_\a h,g_\a P_\a h\rangle = h^* P_\a S_\a P_\a h =
 h^* P_\a [\dot\vp^*\vp'](1,\l_\a) P_\a h
$$
$$
=-h^* P_\a^\bot [\vp^*\vp'](1,\l_\a) P_\a h = -h^* P_\a^\bot [(\vp')^*\vp](1,\l_\a) P_\a h =
0,
$$
since $\vp(1,\l_\a)P_\a=0$. Therefore, $ P_\a h=0$ and $\vp(1,\l_\a)P_\a^\bot h=0$, i.e.
$P_\a^\bot h\in \cE_\a$ and $h=0$.

\no (ii) Note that $\det\vp(1,\l)=0$ if and only if $\vp(1,\l)h=0$ for some $h\in\C^N$, $h\ne
0$, i.e. if and only if $\l$ is an eigenvalue of the operator $H$. Let $\l\!-\!\l_\a=\m$. Due
to $\vp(1,\l_\a)P_\a=0$, we have
$$
\vp(1,\l)= (\m\dot{\vp}(1,\l_\a)+O(\m^2))P_\a+(\vp(1,\l_\a)+O(\m))P_\a^\bot =(Z_\a+O(\m))(\m
P_\a + P_\a^\bot)
$$
as $\m\to 0$. This implies (\ref{Vp-1=}), since $\det Z_\a\ne 0$. Moreover,
$$
\det\vp(1,\l)=\det (Z_\a +O(\m))\det(\m P_\a +P_\a^\bot)=(\det Z_\a +O(\m))\m^{k_\a}\quad {\rm
as }\ \  \m\to 0,
$$
i.e. the multiplicity of the root $\l_\a$ is equal to $k_\a$.

\no (iii) Note that $\vp^{-1}(1,\l)$ is analytic outside the set $\bigcup_{\a\ge 1}\{\l_\a\}$.
Fix some $\a\ge 1$. Due to $\xi(\l_\a)P_\a=0$, we have $\mu^{-1}\xi(\l)P_\a=O(1)$ as $\mu\to
0$. It follows from (\ref{Vp-1=}) that $\xi(\l)\vp^{-1}(1,\l)$ is bounded near $\l_\a$ for
each $\a\ge 1$. Therefore, $\xi(\l)\vp^{-1}(1,\l)$ is entire.
\end{proof}

Recall that we use the notations $\wt{\vp}(x,\l)=\vp(x,\l,\wt{V})$,
$\wt{\cE}_\a=\cE_\a(\wt{V})$ and so on.

\begin{proposition}
\label{*VP1=VP2} (i)  Let $\wt{V}\in \Iso(V)$ for some $V=V^*\in L^1(0,1)$
 and let $\wt{\cE}_\a=\cE_\a$ for all $\a\ge 1$. Then
$$
\wt{\vp}(1,\l)=\vp(1,\l)\quad {for\ all}\quad \l\in\C.
$$

\no (ii) Let, in addition, $\wt{g}_\a=g_\a$ for all $\a\ge 1$. Then
$$
[\wt{\vp}'\wt{\vp}^{-1}](1,\l)=[\vp'\vp^{-1}](1,\l)\quad {for\ all}\quad  \l\!\in\!\C.
$$
\end{proposition}
\begin{proof}
(i) Due to $\wt{P}_\a=P_\a$ and Lemma \nolinebreak \ref{*ZaProp} \nolinebreak (iii), the
function $[\wt{\vp}\vp^{-1}](1,\l)$ is entire. Moreover, asymptotics (\ref{VpAsympt}) gives
$[\wt{\vp}\vp^{-1}](1,z^2)\!=\!I+O(|z|^{-1})$ as $|z|=\pi(n+\frac{1}{2})\to\infty$. Using the
Liouville Theorem, we obtain $[\wt{\vp}\vp^{-1}](1,\l)=I$ for all $\l\in\C$.

\no (ii) Put
$$
f(\l)=[(\wt{\vp}'\wt{\vp}^{-1})\!-\!(\vp'\vp^{-1})](1,\l)=[(\wt{\vp}'\!-\!\vp')\vp^{-1}](1,\l).
$$
Firstly, we prove that the function $f(\l)$ is entire. Due to Lemma \ref{*ZaProp} (iii), it is
sufficient to check that $[\wt{\vp}'-\vp'](1,\l_\a)P_\a=0$ for all $\a\ge 1$. Recall that
$Z_\a=\dot\vp(1,\l_\a)P_\a +\vp(1,\l_\a)P_\a^\bot$ and $\det Z_\a\ne 0$. Using (\ref{LIdent})
and (\ref{VpIdent}), we get
$$
Z_\a^*[\wt{\vp}'-\vp'](1,\l_\a)P_\a = \wt{G}_\a-G_\a +
P_\a^\bot[(\wt{\vp}'-\vp')^*\vp](1,\l_\a)P_\a = 0,
$$
since $\wt{G}_\a=G_\a$ and $\vp(1,\l_\a)P_\a=0$. This gives $[\wt{\vp}'-\vp'](1,\l_\a)P_\a=0$
for all $\a\ge 1$.

Secondly, note that asymptotics (\ref{VpAsympt}), (\ref{Vp'Asympt}) yield $f(z^2)=O(1)$ as
$|z|=\pi(n\!+\!\frac{1}{2})\to\infty$, and $f(z^2)\to 0$ as $z\to i\infty$. Using the
Liouville Theorem, we obtain $f(\l)=0$, $\l\!\in\!\C$.
\end{proof}

\begin{proof}[{\bf Proof of Theorem \ref{*UniqThm}.}\ ]
Recall that $\c(x,\l)$ is the solution of the equation $-\c''+V\c=\l\c$ under the conditions
$\c(1,\l)=0$ and $\c'(1,\l)=I$. Introduce the $2N\!\ts\!2N$ matrix
\[
\label{x1} K(x,\l)= {\matr {\wt\vp} {\wt\c} {\wt{\vp}'} {\wt{\c}'}}(x,\l) {\matr {\vp} {\c}
{\vp'} {\c'}}^{\!-1}(x,\l),\quad x\in [0,1],\ {\l\in\C}.
\]
Using identities (\ref{VpIdent}), (\ref{VpCIdent}), we obtain
\[
\label{x2} {\matr {\vp} {\c} {\vp'} {\c'}}^{\!-1}\!\!(x,\l) = {\matr {\vp^{-1}(1,\l)\!\!\!} 0
0 {(\vp^{-1})^*(1,\ol{\l})}} {\matr {\ (\c')^*} {-\c^*} {-(\vp')^*} {\ \vp^*}}(x,\ol{\l}).
\]
Therefore, $K(x,\l)$ satisfies the differential equation
\[
\label{DiffEqK} K'(x,\l)= {\matr 0 {I} {\wt{V}(x)-\l} 0}\cdot K(x,\l) - K(x,\l)\cdot {\matr 0
{I} {V(x)-\l} 0},\quad x\in [0,1].
\]
It follows from Proposition \ref{*VP1=VP2} that $K(x,\l)$ satisfies the initial condition
$$
K(1,\l)= {\matr {\wt{\vp}(1,\l)} 0 {\wt{\vp}'(1,\l)} {I}} {\matr {\ \vp^{-1}(1,\l)} 0
{-[\vp'\vp^{-1}](1,\l)} {I}} = {\matr {I} 0 0 {I}}.
$$
Then, each function $K(x,\cdot)$, $x\in [0,1]$, is entire as a solution of equation
(\ref{DiffEqK}).

Fix some $x\in [0,1]$. Substituting asymptotics (\ref{VpAsympt}), (\ref{Vp'Asympt}) into
(\ref{x1}), (\ref{x2}), we obtain
$$
K(x,z^2)= {\matr {I+O(|z|^{-1})} {O(|z|^{-2})} {O(1)} {I+O(|z|^{-1})}}\quad {\rm as}\ \
|z|=\pi(n+{\textstyle\frac{1}{2}})\to\infty.
$$
and $K(x,z^2)\to I_{2N}$ as $z\to i\infty$. Hence, $K(x,\l)=I_{2N}$ for each $(x,\l)\in
[0,1]\ts\C$. In particular, this gives $\wt{\vp}(x,\l)=\vp(x,\l)$ for all $(x,\l)\in
[0,1]\ts\C$, i.e. $\wt{V}=V$.
\end{proof}


Introduce the subspaces
\[
\label{cEsDef} \cE^\sharp_\a=\Ker \vp^*(1,\l_\a,V),\quad \a\ge 1.
\]
and let $P_\a^\sharp:\C^N\to \cE_\a^\sharp$ be the orthogonal projector. Using
(\ref{VpCIdent}) and (\ref{ChiIdentity}), we obtain $\cE_\a^\sharp(V) =\cE_\a(V^\sharp)$,
$P_\a^\sharp(V)=P_\a(V^\sharp)$ for all $\a\ge 1$, where $V^\sharp(t)=V (1-t)$, $t\in [0,1]$.
\begin{lemma}
\label{*xIdent} Let $V=V^*\in L^1(0,1)$ and $\a\ge 1$. The following identities are fulfilled:
\[
\label{xIdent1} P_\a^\sharp\vp'(1,\l_\a)P_\a=\vp'(1,\l_\a)P_\a,\qquad
\c'(0,\l_\a)\vp'(1,\l_\a)P_\a=P_\a,
\]
\[
\label{xIdent2} \res_{\l=\l_\a}\c^{-1}(0,\l_\a)\cdot \dot{\c}(0,\l_\a)
P_\a^\sharp=P_\a^\sharp.
\]
\end{lemma}
\begin{proof}
Due to (\ref{VpIdent}), we have $[\vp^*\vp'](1,\l_\a)P_\a=[(\vp')^*\vp] (1,\l_\a)P_\a=0$. This
yields the first identity in (\ref{xIdent1}). Let
$\e(x)=\c(x,\l_\a)\vp'(1,\l_\a)P_\a-\vp(x,\l_\a)P_\a$. Then the function $\e$ satisfies the
equation $-\e''+V\e=\l_\a\e$ and
$$
\e(1)=-\vp(1,\l_\a)P_\a=0,\qquad \e'(1)=\vp'(1,\l_\a)P_\a-\vp'(1,\l_\a)P_\a=0.
$$
Therefore, $\e(x)=0$ for all $x\in[0,1]$. Using $\e'(0)=0$, we obtain the second identity in
(\ref{xIdent1}). Furthermore, asymptotics (\ref{Vp-1=}) yields
$$
\res_{\l=\l_\a}\vp^{-1}(1,\l)\cdot \dot{\vp}(1,\l_\a)P_\a = P_\a Z_\a^{-1}\cdot
\dot{\vp}(1,\l_\a)P_\a = P_\a Z_\a^{-1}\cdot Z_\a P_\a =P_\a,
$$
since $Z_\a=\dot{\vp}(1,\l_\a)P_\a+\vp(1,\l_\a)P_\a^\bot$. Applying
this formula with the
potential $V^\sharp$ instead of $V$ and using (\ref{ChiIdentity}) and
$\cP_\a(V^\sharp)=\cP_\a^\sharp(V)$, we obtain (\ref{xIdent2}).
\end{proof}

\begin{proof}[{\bf Proof of Proposition \ref{*ResM=}.}]
Identity (\ref{VpIdent}) gives $m(\l)=m^*(\ol{\l})$. This implies
$$
B_\a= - \res_{\l=\l_\a}m(\l)=B_\a^*
$$
for all $\a\ge 1$. Due to (\ref{Vp-1=}), the function $\c^{-1}(1,\l)=-(\vp^*)^{-1}(1,\ol{\l})$
has a simple pole at each point $\l=\l_\a$. Therefore, the function $m(\l)$ has a simple pole
at $\l=\l_\a$ and
$$
B_\a=\c'(0,\l_\a)\cdot \res_{\l=\l_\a} (\vp^*)^{-1}(1,\ol{\l}) =
\c'(0,\l_\a)
(Z_\a^*)^{-1}P_\a\,.
$$
This yields $B_\a\big|_{\C^N\ominus\cE_\a}=0$. Recall that $G_\a=P_\a
S_\a P_\a=P_\a
[\dot{\vp}^*\vp'](1,\l_\a)P_\a$ (see (\ref{LIdent})). Hence,
$$
B_\a G_\a = B_\a [\dot{\vp}^*\vp'](1,\l_\a)P_\a = \c'(0,\l_\a)\cdot
\res_{\l=\l_\a}
\c^{-1}(0,\l)\cdot \dot{\c}(0,\l_\a)\vp'(1,\l_\a)P_\a,
$$
where we have used the identity $\dot{\vp}^*(1,\l_\a)= -\dot\c(0,\l_\a)$. Applying Lemma
\ref{*xIdent}, we obtain
$$
B_\a G_\a = \c'(0,\l_\a) \cdot \res_{\l=\l_\a} \c^{-1}(0,\l)\cdot
\dot{\c}(0,\l_\a)P_\a^\sharp\vp'(1,\l_\a)P_\a
$$
$$
= \c'(0,\l_\a)P_\a^\sharp\vp'(1,\l_\a)P_\a = \c'(0,\l_\a)\vp'(1,\l_\a)
P_\a = P_\a,\quad \a\ge
1.
$$
Therefore, $B_\a\big|_{\cE_\a}=g_\a^{-1}$ for all $\a\ge 1$.
\end{proof}

\section{Isospectral transforms}
\setcounter{equation}{0}

Let $V=V^*\in L^1(0,1)$ and $\a\ge 1$. We will use the notations
$$
\vp_\a(x)=\vp(x,\l_\a),\qquad \ S_\a(x)=\int_0^x [\vp_\a^*\vp_\a](t)dt,\qquad S_\a=S_\a(1).
$$

Recall that the residues of the Weyl-Titchmarsh function are given by
\[
\label{Bdef} B_\a= -\res_{\l=\l_\a}m(\l,V)=B_\a^*,\qquad B_\a\Big|_{\cE_\a}=g_\a^{-1},\quad
B_\a\Big|_{\C^N\ominus\cE_\a}=0
\]
(see Proposition \nolinebreak \ref{*ResM=}). Note that $B_\a\ge 0$  and {\bf the matrix $\bf
B_\a$ is uniquely determined by the pair $\bf \{\cE_\a\,, g_\a\}$ and visa versa.} In
particular, $\cE_\a=\C^N\ominus\Ker B_\a$.
\begin{lemma}
Let $V=V^*\in L^1(0,1)$, $\a\ge 1$ and let $e_\a$ be some $k_\a\!\ts\!N$ matrix whose columns
form the basis of the subspace $\cE_\a$. Then
\[
\label{Ba=} B_\a= e_\a\left[e_\a^*S_\a e_\a\right]^{-1}e_\a^*.
\]
\end{lemma}
\begin{proof}
If the columns of $e_\a$ form the orthonormal basis in $\cE_\a$, then (\ref{Bdef}) and
(\ref{Ba=}) are equivalent. Let $e_\a'$ be another $k_\a\!\ts\!N$ matrix whose columns form
the basis of $\cE_\a$. Then, $e_\a=e'_\a U$ for some $k_\a\!\ts\!k_\a$ matrix $U$ such that
$\det U\ne 0$. The simple calculation shows that
$$
B_\a(V)= e'_\a U \left[U^*(e'_\a)^*S_\a(1) e'_\a U\right]^{-1} U^*(e'_\a)^* =
e'_\a\left[(e'_\a)^*S_\a(1) e'_\a\right]^{-1}(e'_\a)^*.
$$
Hence, (\ref{Ba=}) doesn't depend on the choice of $e_\a$.
\end{proof}

Introduce the matrices
$$
D_\a =S_\a^{-1}-B_\a,\quad \a\ge 1.
$$
\begin{lemma}
\label{*DaLemma} Let $V=V^*\in L^1(0,1)$, $\a\ge 1$ and the subspace $\cF_\a\ss\C^N$ be given
by (\ref{FaDef}). Then $D_\a=D_\a^*\ge 0$, $\cF_\a=\C^N\ominus\Ker D_\a$ and $\dim
\cF_\a=N-k_\a$. Moreover, $\cF_\a\cap\cE_\a=\{0\}$.
\end{lemma}
\begin{proof}
Without loss of generality, we can assume that
$$
e_\a= {\vect {I_{k_\a}} 0};\qquad S_\a = {\matr s p {p^*} q},\qquad B_\a=e_\a (e_\a^* S_\a
e_\a)^{-1} e_\a^* = {\matr {s^{-1}} 0 0 0},
$$
where $s=s^*$ is a $k_\a\!\ts\!k_\a$ matrix, $p$ is $(N\!-\!k_\a)\!\ts\!k_\a$ and $q=q^*$ is
$(N\!-\!k_\a)\!\ts\!(N\!-\!k_\a)$. Note that $S_\a>0$ yields $s>0$ and $q>0$. Due to the
Frobenius formula for the inverse matrix (see \cite{G}, Ch. 2.5), we have
$$
S_\a^{-1}= {\matr {(s-pq^{-1}p^*)^{-1}} {-s^{-1}p(q-p^*s^{-1}p)^{-1}}
{-q^{-1}p^*(s-pq^{-1}p^*)^{-1}} {(q-p^*s^{-1}p)^{-1}}}.
$$
Note that $q^{-1}p^*(s-pq^{-1}p^*)^{-1}=(q-p^*s^{-1}p)^{-1}p^*s^{-1}$, since
$S_\a^{-1}=(S_\a^{-1})^*$. Together with the identity
$(s-pq^{-1}p^*)^{-1}-s^{-1}=s^{-1}pq^{-1}p^*(s-pq^{-1}p^*)^{-1}$, this yields
$$
D_\a=S_\a^{-1}-B_\a = {\matr {s^{-1}pq^{-1}p^*(s-pq^{-1}p^*)^{-1}}
{-s^{-1}p(q-p^*s^{-1}p)^{-1}} {-q^{-1}p^*(s-pq^{-1}p^*)^{-1}} {(q-p^*s^{-1}p)^{-1}}}
$$
$$
\ \ \ \ \  = {\vect {-s^{-1}p} {I_{N-k_\a}}}(q-p^*s^{-1}p)^{-1}\lt(-p^*s^{-1}\ \
I_{N-k_\a}\rt).
$$
This implies $D_\a\ge 0$, $\rank D_\a=N-k_\a$ and $\dim\Ker D_\a=k_\a$. Moreover, the identity
$$
\lt(-p^*s^{-1}\ \ I_{N-k_\a}\rt){\vect s {p^*}}=0
$$
yields $S_\a(\cE_\a)\ss\Ker D_\a$. Recall that $S_\a=S_\a^*>0$. Using
$\dim\cE_\a=k_\a=\dim\Ker D_\a$, we deduce that $S_\a(\cE_\a)= \Ker D_\a$ and so
$\cF_\a=\C^N\ominus \Ker D_\a$. Since $D_\a+B_\a>0$, we have $\Ker D_\a \cap \Ker B_\a =
(\C^N\ominus\cF_\a)\cap(\C^N\ominus\cE_\a)=\{0\}$. Together with $\dim\cF_\a + \dim\cE_\a =
N$, this implies $\cF_\a\cap\cE_\a=\{0\}$.
\end{proof}
\begin{corollary}
\label{*FaCor} Let $V=V^*\in L^1(0,1)$. Then $\cF_\a=\C^N\ominus
[\dot{\vp}^*(1,\l_\a)](\cE_\a^\sharp)$ for all $\a\ge 1$, where the subspace
$\cE_\a^\sharp\ss\C^N$ is given by (\ref{cEsDef}).
\end{corollary}
\begin{proof}
Recall that $\cF_\a=\C^N\ominus S_\a(\cE_\a)$ and $S_\a=\int_0^1[\vp^*\vp](t,\l_\a)dt=
[\dot{\vp}^*\vp'-(\dot{\vp}')^*\vp](1,\l_\a)$. Since $\vp(1,\l_\a)\big|_{\,\cE_\a}=0$, we
obtain
$$
S_\a(\cE_\a)=\left[\dot{\vp}^*(1,\l_\a)\vp'(1,\l_\a)\right](\cE_\a).
$$
Using (\ref{xIdent1}), we deduce that $[\vp'(1,\l_\a)](\cE_\a)=\cE_\a^\sharp$. Hence,
$S_\a(\cE_\a)=[\dot{\vp}^*(1,\l_\a)](\cS_\a^\sharp)$.
\end{proof}

The following Theorem gives the explicit formula for the isospectral transform of the
potential $V$, changing only the matrix $B_\a(V)$.

\begin{theorem}
\label{*TransfThm} Let the potential $V\!=\!V^*\in L^1(0,1)$ and the matrix $B=B^*$ be such
that
\[
\label{Acond} B\ge 0,\quad \rank B = k_\a \quad {and}\quad \cE^{(B)}\cap \cF_\a=\{0\},\ \
{where}\ \ \cE^{(B)}=\C^N\ominus \Ker B,
\]
for some $\a\ge 1$. Denote
\[
\label{WtVDef} \wt V(x)=V(x)-2[\vp_\a K \vp^*_\a]'(x) = \wt{V}^{*}(x),\quad x\in [0,1],
\]
where
\[
\label{KandADef}
K(x)=A(I\!+\!S_\a(x)A)^{-1}\!=
K^*(x)\quad {and} \quad A=B - B_\a.
\]
Then $(\wt V - V)'\in L^1(0,1)$, $\wt{B}=B_\b$ for all $\b\ne\a$, and
$$
\wt{B}_\a=B,\qquad \wt{\cE}_\a=\cE^{(B)},\qquad \wt{\cF}_\a=\cF_\a.
$$
\end{theorem}

In order to prove Theorem \ref{*TransfThm}, we need two preliminary Lemmas.

\begin{lemma}
\label{*VpKVp=0} Let $V\!=\!V^*\in L^1(0,1)$, $\a\!\ge\! 1$ and $B\!=\!B^*$ be such that
(\ref{Acond}) hold true. Then

\no (i) $\det (I\!+\!S_\a(x)A)\ne 0$ for all $x\in [0,1]$, where $A=B-B_\a$.

\no (ii) the identity $[\vp_\a K \vp_\a^*](1)=0$ is fulfilled, where $K$ is given by
(\ref{KandADef}).
\end{lemma}
\begin{proof}
(i) Note that $I+S_\a(x)A=S_\a(x)(S_\a^{-1}(x)+A)$. For all $x\in [0,1]$ we have
$$
S_\a^{-1}(x)+A> S_\a^{-1}(1)+A = (S_\a^{-1}(1)\!-\!B_\a)+ (B_\a\!+\!A)= D_\a + B \ge 0,
$$
since $D_\a\ge 0$ and $B\ge 0$. Lemma \ref{*DaLemma} and (\ref{Acond}) give $\rank D_\a =
N-k_\a$, $\rank B = k_\a$ and $(\C^N\ominus \Ker D_\a)\cap (\C^N\ominus \Ker B)=\{0\}$. Hence,
$\Ker D_\a \cap \Ker B = \{0\}$ and $D_\a + B>0$.

\no (ii)  Let the columns of the matrix $e_\a$ form some orthonormal basis in $E_\a$\,.
Without loss of generality, we can assume that
$$
e_\a= {\vect {I_{k_\a}} 0};\quad S_\a(1) = {\matr s p {p^*} q},\quad B= {\matr b c {c^*} d},
$$
where $s=s^*$ and $b=b^*$ are $k_\a\!\ts\!k_\a$ matrices, $p$ and $c$ are
$(N\!-\!k_\a)\!\ts\!k_\a$ matrices, $q=q^*$ and $d=d^*$ are $(N\!-\!k_\a)\!\ts\!(N\!-\!k_\a)$
matrices. Then,
$$
B_\a=e_\a (e_\a^* S_\a(1) e_\a)^{-1} e_\a^* = {\matr {s^{-1}} 0 0 0} ,\qquad A = {\matr
{b-s^{-1}} c {c^*} d}.
$$

Firstly, let $b>0$. In this case, the assumption $\rank B=k_\a$ gives $d=c^*b^{-1}c$. Note
that
$$
I+S_\a(1) A = {\matr {sb+pc^*} {sc+pd} {p^*(b-s^{-1})+qc^*} {I_{N-k_\a}+p^*c+qd} }
$$
and
\[
\label{xRow1}
 \lt({sb+pc^*}\ \ {sc+pd} \rt)= (s+pc^*b^{-1})\cdot \lt(b\ \ c\rt).
\]
Moreover,
\[
\label{xRow2} \lt(c^*\ \ d\rt)=c^*b^{-1}\cdot \lt(b\ \ c\rt).
\]
Due to (\ref{xRow1}), (\ref{xRow2}), for each index $j\in[k_\a\!+\!1,N]$, the first $k_\a$
rows of the matrix $I+S_\a(1)A$ and the $j$-th row of $A$ are linearly dependent. Using the
Cramer formula\begin{footnote}{ Let $X=[x_{j,m}]_{j,m=1}^N$, $Y=[y_{j,m}]_{j,m=1}^N$ be
$N\!\ts\!N$ matrices and let $\det Y\ne 0$, $Z=XY^{-1}=[z_{j,m}]_{j,m=1}^N$. Then
$$
z_{j,m}=\frac{1}{\det Y}\left|\begin{array}{ccccc} y_{1,1} & y_{1,2} & \dots & y_{1,N-1} &
y_{1,N} \cr \dots &  \dots &  \dots &  \dots & \dots \cr y_{m-1,1} & y_{m-1,2} & \dots &
y_{m-1,N-1} & y_{m-1,N} \cr x_{j,1} & x_{j,2} & \dots & x_{j,N-1} & x_{j,N} \cr y_{m+1,1} &
y_{m+1,2} & \dots & y_{m+1,N-1} & y_{m+1,N} \cr \dots &  \dots &  \dots &  \dots & \dots \cr
y_{N,1} & y_{N,2} & \dots & y_{N,N-1} & y_{N,N}
\end{array}\right|\qquad {\rm for\ all}\quad j,m\in[1,N]\quad {\rm (see\ \cite{G},\ Ch.\ 1.3).}
$$
}\end{footnote}, we deduce that the matrix $K(1)=A(I+S_\a(1) A)^{-1}$ has the form
$$
K(1)= {\matr {...} {...} {c^*} d}{\matr {sb+pc^*} {sc+pd} {...} {...}}^{-1}= {\matr {...}
{...} {...} 0}.
$$
Recall that $\vp_\a(1)e_\a=0$. Hence,  $\vp_\a(1)= {\matr 0 {...} 0 {...}}$. This yields
$[\vp_\a^* K\vp_\a](1) =0$.

Secondly, let $\det b=0$. Since $K(1)$ is a continuous function of  $B$, we deduce that it has
the same form as before and $[\vp_\a^*K\vp_\a](1)=0$.
\end{proof}

Introduce the matrices
$$
S(x,\l)=\int_0^x [\vp^*\vp](t,\l)dt,\qquad \wt{S}(x,\l)=\int_0^x[\wt{\vp}^*\wt{\vp}](t,\l)dt.
$$
\begin{lemma}
\label{*TildeLemma} Let the potential $V\!=\!V^*\in L^1(0,1)$, the number $\a\!\ge\!1$ and the
matrix $B=B^*$ be such that conditions (\ref{Acond}) are fulfilled. Then
$$
\wt{\vp}(x,\l)=\vp(x,\l)-[\vp_\a KT](x,\l)\quad {and}\quad \wt{S}(x,\l)=
S(x,\l)-[T^*KT](x,\l),
$$
where
$$
T(x,\l)=\int_0^x\vp_\a^*(t)\vp(t,\l)dt
$$
and the potential $\wt{V}$ is given by (\ref{WtVDef}).
\end{lemma}
\begin{proof}
Let $\e=\vp-\vp_\a KT$ (we omit $x$ and $\l$ for short). Using the identity
$K'=-K\vp_\a^*\vp_\a K$, we obtain $(KT)'= K\vp_\a^*\e$. Therefore,
$$
\e'= \vp' - \vp'_\a KT - \vp_\a K\vp_\a^*\e,
$$
$$
\e''=\vp''-\vp''_\a KT - \vp'_\a K\vp_\a^*\e - (\vp_\a K\vp_\a^*)'\e - \vp_\a K\vp_\a^*
\left(\vp' - \vp'_\a KT - \vp_\a K\vp_\a^*\e\right).
$$
Recall that $\vp_\a^*\vp_\a'=(\vp_\a')^*\vp_\a$  (see (\ref{VpIdent})). Hence,
$$
\e''= \vp''-\vp''_\a KT - 2(\vp_\a K\vp_\a^*)'\e + \vp_\a K ((\vp'_\a)^*\vp-\vp_\a^*\vp').
$$
Note that $\vp''=(V\!-\!\l)\vp$, $\vp''_\a=(V\!-\!\l_\a)\vp_\a$ and
$(\vp'_\a)^*\vp-\vp_\a^*\vp'=(\l\!-\!\l_\a)T$. This gives
$$
\e''= (V\!-\!\l)\e-2(\vp_\a K\vp_\a^*)'\e = (\wt{V}\!-\!\l)\e.
$$
Since $\e(0)=0$ and $\e'(0)=\vp'(0)=I$, we deduce $\e=\wt{\vp}$. Furthermore,
$$
\e^*\e=\vp^*\vp -\vp^*\vp_\a KT + T^*K\vp_\a^*\vp_\a KT - T^*K\vp_\a^*\vp = \vp^*\vp -
[T^*KT]'.
$$
This yields $\wt S = S - T^*KT$.
\end{proof}

\begin{proof}[{\bf Proof of Theorem \ref{*TransfThm}}.\ ]
Note that $(\vp_\a K\vp_\a^*)''\in L^1(0,1)$, since $\vp_\a'',K''\in L^1(0,1)$. Fix some
$\b\ne\a$. Let $\vp_\b(x)=\vp(x,\l_\b,V)$, $\wt{\vp}_\b(x)=\vp(x,\l_\b,\wt{V})$ and so on.
Using Lemma \ref{*TildeLemma} and the identity
\[
\label{Tb=} T_\b(1)=\int_0^1[\vp_\a^*\vp_\b](t)dt=
\frac{[\vp_\a^*\vp'_\b-(\vp'_\a)^*\vp_\b](1)}{\l_\a-\l_\b},
\]
we obtain
\[
\label{WtVpB} \wt{\vp}_\b(1)= \vp_\b(1)-\frac{[\vp_\a K
(\vp_\a^*\vp'_\b-(\vp'_\a)^*\vp_\b)](1)}{\l_\a-\l_\b} = \lt(I + \frac{[\vp_\a
K(\vp'_\a)^*](1)}{\l_\a-\l_\b}\rt)\vp_\b(1),
\]
where we have used $[\vp_\a K \vp_\a^*](1)=0$ (see Lemma \ref{*VpKVp=0} (ii)). Therefore,
$\l_\b$ is a root of the (scalar) entire function $\wt{w}(\l)=\det\wt{\vp}(1,\l)$ of the
multiplicity at least $k_\b$. Furthermore,
\[
\label{WtVpA} \wt{\vp}_\a(1)= \vp_\a(1)-[\vp_\a K S_\a](1)=\vp_\a(1)\left(I-[KS_\a](1)\right).
\]
Hence, $\l_\a$ is a root of $\wt{w}(\l)$ of the multiplicity at least $k_\a$. Using Lemma
\ref{*ZaProp} (ii), we deduce that $[\wt{w}w^{-1}](\l)$ is an entire function, where
$w(\l)=\det\vp(1,\l)$. Note that $[\wt{w}w^{-1}](z^2)=1+O(|z|^{-1})$ as
$|z|=\pi^2(n\!+\!\frac{1}{2})^2\to\infty$ (see (\ref{VpAsympt})). Hence, $\wt{w}(\l)=w(\l)$,
$\l\in\C$. In other words, each $\l_\b$\,, $\b\ge 1$, is an eigenvalue of the operator
$\wt{H}\p=-\p''+\wt{V}\p$ of the multiplicity $k_\b$ and there are no other eigenvalues.

Let $\b\ne\a$. Identity (\ref{WtVpB}) gives $\wt{\cE}_\b\supset\cE_\b$. Since $\dim
\wt{\cE}_\b=k_\b=\dim\cE_\b$, we obtain $\wt{\cE}_\b=\cE_\b$\,. Let $e_\b$ be some
$k_\b\!\ts\!N$ matrix whose columns form the basis of $E_\b$. Recall that
$\wt{B}_\b=e_\b(e_\b^*\wt{S}_\b e_\b)^{-1}e_\b^*$. Using Lemma \ref{*TildeLemma} and
(\ref{Tb=}), we obtain
$$
\wt{S}_\b = S_\b -
\frac{[((\vp'_\b)^*\vp_\a-\vp_\b^*\vp'_\a)K(\vp_\a^*\vp'_\b-(\vp'_\a)^*\vp_\b)](1)}
{(\l_\a-\l_\b)^2}.
$$
Note that $\vp_\b(1)e_\b=0$, since the columns of $e_\b$ belong to $\cE_\b=\Ker \vp_\b(1)$.
Due to Lemma \nolinebreak \ref{*VpKVp=0}, $[\vp_\a^* K \vp_\a](1)=0$. Therefore,
$e_\b^*\wt{S}_\b e_\b = e_\b^*{S}_\b e_\b$ and
$$
\wt{B}_\b=e_\b(e_\b^*\wt{S}_\b e_\b)^{-1}e_\b^*= e_\b(e_\b^* S_\b e_\b)^{-1}e_\b^* = B_\b,
\quad \b\ne\a.
$$

We will show that $\wt{B}_\a=B$. Let the columns of the matrix $e_\a$ form some orthonormal
basis in $E_\a$. Without loss of generality, we can assume that
$$
e_\a= {\vect {I_{k_\a}} 0};\quad S_\a = {\matr s p {p^*} q},\quad B= {\matr b c {c^*} d},
$$
where $s=s^*$ and $b=b^*$ are $k_\a\!\ts\!k_\a$ matrices, $p$ and $c$ are
$(N\!-\!k_\a)\!\ts\!k_\a$ matrices, $q=q^*$ and $d=d^*$ are $(N\!-\!k_\a)\!\ts\!(N\!-\!k_\a)$
matrices. Note that
$$
B_\a=e_\a (e_\a^* S e_\a)^{-1} e_\a^* = {\matr {s^{-1}} 0 0 0} ,\qquad A = {\matr {b-s^{-1}} c
{c^*} d}.
$$

Firstly, let $b>0$. In this case, $\rank B=k_\a$ implies $d=c^*b^{-1}c$. It follows from
identity (\ref{WtVpA}) that the columns of the matrix
$$
\wt{e}_\a=(I-KS_\a)^{-1}e_\a = (I-A(S_\a^{-1}\!+\!A)^{-1})^{-1}e_\a = (I+AS_\a)e_\a
$$
form a basis in $\wt{E}_\a$ (recall that $\det (I\!+\!AS_\a)=\det (I\!+\!S_\a A)\!\ne\!0$ due
to Lemma \ref{*VpKVp=0} (i)). We have
$$
\wt{e}_\a = {\vect {I_{k_\a}+(b-s^{-1})s+cp^*} {c^*s+dp^*}} = {\vect {bs+cp^*}
{c^*s+c^*b^{-1}cp^*} } = {\vect {I_{k_\a}} {c^*b^{-1}}}(bs+cp^*).
$$
Furthermore, Lemma \ref{*TildeLemma} gives
$$
\wt{S}_\a=S_\a-S_\a^*KS_\a = S_\a (I-KS_\a)=S_\a(I+AS_\a)^{-1}.
$$
Hence,
$$
{\wt{e}_\a}^{\,\,*}\wt{S}_\a\wt{e}_\a = {\wt{e}_\a}^{\,\,*} S_\a e_\a = (sb+pc^*)\lt(I_{k_\a}\
\ b^{-1}c\rt){\vect s {p^*}}= (sb+pc^*)b^{-1}(bs+cp^*).
$$
We obtain
$$
\wt{B}_\a=\wt{e}_\a ({\wt{e}_\a}^{\,\,*}\wt{S}_\a\wt{e}_\a)^{-1}{\wt{e}_\a}^{\,\,*} = {\vect
{I_{k_\a}} {c^*b^{-1}}}\cdot b\cdot \lt(I_{k_\a}\ \ b^{-1}c\rt) = {\matr b c {c^*}
{c^*b^{-1}c}} = B.
$$

Secondly, let $\det b=0$. Note that $\wt{e}_\a$ and $\wt{S}_\a$ are continuous functions of
$B$. Therefore, $\wt{B}_\a=B$ due to the arguments given above and the continuity of
$\wt{B}_\a$ as a function of $B$. Since $\wt{B}_\a=B$, we obtain $\wt{\cE}_\a=\C^N\ominus\Ker
B = \cE^{(B)}$. Note that
$$
\wt{S}_\a^{-1}=\left(S_\a(I+AS_\a)^{-1}\right)^{-1}=S_\a^{-1}+A.
$$
This implies
$$
\wt{D}_\a=\wt{S}_\a^{-1}-\wt{B}_\a=S_\a^{-1}+A-B=S_\a^{-1}-B_\a=D_\a.
$$
In particular, we have $\wt{\cF}_\a=\C^N\ominus\Ker \wt{D}_\a=\C^N\ominus \Ker D_\a=\cF_\a$.
\end{proof}

\begin{proof}[{\bf Proof of Theorem \ref{*EGTransThm}.\ }]
Due to Theorem \ref{*UniqThm}, the mapping $\F_\a$ is one-to-one. We prove that $\F_\a$ is
onto. Let the subspace $\cE\ss\C^N$ be such that $\dim\cE=k_\a$, $\cE\cap\cF_\a=\{0\}$ and let
$g=g^*>0$ be the operator in $\cE$. We define the matrix $B=B^*\ge 0$ by
$$
B\Big|_{\cE}=g^{-1},\qquad B\Big|_{\C^N\ominus \cE}=0.
$$
Since $\cE^{(B)}=\C^N\ominus\Ker B = {\cE}$, conditions (\ref{Acond}) are fulfilled. Let
$\wt{V}$ be given by formula (\ref{WtVDef}). It follows from Theorem \ref{*TransfThm} that
$\wt{V}\in\Iso(V)$, $\wt{B}_\b=B_\b$ for all $\b\ne\a$ and $\wt{B}_\a=B$. In view of
definition (\ref{Bdef}), this yields $\wt{\cE}_\b=\cE_\b$, $\wt{g}_\b=g_\b$ for all $\b\ne\a$
and $\wt{\cE}_\a={\cE}$, $\wt{g}_\a=g$.
\end{proof}

\begin{proof}[{\bf Proof of Proposition \ref{*FaProp}.\ }]
Since $\dim \cF_\a=\dim \wt{\cF}_\a = N\!-\!k_\a$, we can fix some matrix $B=B^*\ge 0$ such
that $\rank B=k_\a$, $\cE^{(B)}\cap \cF_\a=\{0\}$ and $\cE^{(B)}\cap\wt{\cF}_\a=\{0\}$, where
$\cE^{(B)}=\C^N\ominus \Ker B$. Using Theorem \nolinebreak \ref{*TransfThm}, we construct
potentials $V_B,\wt{V}_B\in\Iso(V)$ such that
$$
B_\b(V_B)=B_\b,\ \ B_\b(\wt{V}_B)=\wt{B}_\b,\ \ \b\ne\a,\qquad B_\a(V_B)=B=B_\a(\wt{V}_B).
$$
Note that $\cF_\a(V_B)=\cF_\a$ and $\cF_\a(\wt{V}_B)=\wt{\cF}_\a$. It is clear that
$$
\cE_\b(V_B)=\cE_\b=\wt{\cE}_\b=\cE_\b(\wt{V}_B),\ \ \b\ne\a,\quad {\rm and}\quad
\cE_\a(V_B)=\cE^{(B)}=\cE_\a(\wt{V}_B).
$$
Due to Proposition \ref{*VP1=VP2} (i), we have $\vp(1,\l,V_B)=\vp(1,\l,\wt{V}_B)$ for all
$\l\in\C$. In particular, $\dot{\vp}^*(1,\l_\a,V_B)=\dot{\vp}^*(1,\l_\a,\wt{V}_B)$ and
$\cE_\a^\sharp(V_B)=\Ker \vp^*(1,\l_\a,V_B)= \Ker
\vp^*(1,\l_\a,\wt{V}_B)=\cE_\a^\sharp(\wt{V}_B)$. Corollary \ref{*FaCor} implies
$\cF_\a(V_B)=\cF_\a(\wt{V}_B)$. Therefore, $\cF_\a=\wt{\cF}_\a$. Together with Lemma
\ref{*DaLemma} this gives $\wt{\cE}_\a\cap\cF_\a=\wt{\cE}_\a\cap\wt{\cF}_\a=\{0\}$.
\end{proof}

\begin{proof}[{\bf Proof of Proposition \ref{*ExProp}}]
Let $\cE_1\cap\cE_2\ne\{0\}$. Due to $\dim\cE_1=\dim\cE_2=1$, this is equivalent to
$\cE_1=\cE_2$. Fix some vector $h\in \cE_1=\cE_2$ such that $\|h\|=1$. Let
$$
w(\l)=\det \vp(1,\l),\qquad u(\l)= [h^* \vp(1,\l) h]^2,\quad \l\in\C.
$$
It follows from (\ref{VpAsympt}) that both functions $w$, $u$ have the same asymptotics
$$
w(z^2),u(z^2)= \frac{\sin^2 z}{z^2}\cdot(1+O(|z|^{-1})) \quad {\rm as}\
|z|=\pi(n\!+\!{\textstyle\frac{1}{2}})\to\infty.
$$
Therefore, by Rouch\'e's Theorem, $w(\l)$ and $u(\l)$ have the same number of roots counted
with multiplicity in the disc $\{|\l|<\pi^2(N\!+\!\frac{1}{2})^2\}$, if $N$ is sufficiently
large. Due to Lemma \nolinebreak \ref{*ZaProp} \nolinebreak (ii), $w(\l)$ has simple roots at
the points $\l_1$, $\l_2$ and $w(\l)$ has double roots at the points $\l_\a$, $\a\ge 3$.
Moreover, $w(\l)$ has no other roots. On the other hand, each point $\l_\a$, $\a\ge 1$, is a
double root of $u(\l)$. This is a contradiction. Hence, $\cE_1\cap\cE_2= \{0\}$.

Suppose that $\cF_1\ne\cE_2$, i.e. $\cF_1\cap\cE_2=\{0\}$. Then, using Theorem
\ref{*EGTransThm}, we can construct the potential $\wt{V}\in\Iso(V)$ such that
$\wt{\cE}_1=\cE_2=\wt{\cE}_2$. Due to the arguments given above, this is impossible.
Therefore, $\cF_1=\cE_2$. The proof of $\cF_2=\cE_1$ is similar.
\end{proof}

\no {\bf Acknowledgments.} Dmitry Chelkak was partly supported by grants VNP Minobrazovaniya
3.1--\,4733 and "Vedushie nauchnye shkoly". Evgeny Korotyaev was partly supported by DFG
project BR691/23-1. Some part of this paper was written at the Mittag-Leffler Institute,
Stockholm. The authors are grateful to the Institute for the hospitality. The authors would
like to thank Markus Klein for useful discussions.


\begin{thebibliography}{CHGL}
\bibitem [Ca] {Ca} Carlson R.: An inverse problem for the matrix
Schr\"odinger equation. J. Math. Anal. Appl. 267 (2002), no. 2, 564--575.
\bibitem [CK] {CK} Chelkak D., Korotyaev E.:
Spectral estimates for Schr\"odinger operators with periodic matrix potentials on the real
line. Preprint Institut Mittag-Leffler No. 06, 2005/2006 fall. To be published in
International Mathematics Research Noticies.
\bibitem [Ch] {Ch} Chern Hua-Huai: On the construction of isospectral vectorial
Sturm-Liouville differential equation. Preprint 1998.
\bibitem [ChSh] {ChSh} Chern Hua-Huai; Shen Chao-Liang:
On the $n$-dimensional Ambarzumyan's theorem. Inverse Problems {\bf 13}(1), 15--18 (1997).
\bibitem [CHGL] {CHGL} Clark S., Holden H., Gesztesy, F., Levitan, B.:
Borg-type theorem for matrix-valued Schr\"odinger operators. J. Diff. Eqs. {\bf 167}, 181--210
(2000).
\bibitem [G] {G} Gantmacher F. R.: {\it The theory of matrices. Vol. 1.}
Translated from the Russian by K. A. Hirsch. Reprint of the 1959 translation.
AMS Chelsea Publishing, Providence, RI, 1998.
\bibitem [JL1] {JL1} Jodeit M. Jr., Levitan B. M.:
Isospectral vector-valued Sturm-Liouville problems. Lett. Math. Phys. {\bf 43}(2), 117--122
(1998).
\bibitem [JL2] {JL2} Jodeit M. Jr., Levitan B. M.:
A characterization of some even vector-valued Sturm-Liouville problems. Mat. Fiz. Anal. Geom.
{\bf 5}(3-4), 166--181 (1998).
\bibitem [L] {L} Levitan B.: {\it Inverse Sturm-Liouville problems}.
Utrecht: VNU Science Press, 1987.
\bibitem [Mal] {Mal} Malamud M. M.:
Uniqueness of the matrix Sturm-Liouville equation given a part of the monodromy matrix, and
Borg type results. {\it Sturm-Liouville Theory: Past and Present}, 237--270. Birkhauser,
Basel, 2005.
\bibitem [M] {M} Marchenko V.: {\it Sturm-Liouville operator and
applications}. Basel: Birkh\"auser, 1986.
\bibitem [PT] {PT} P\"oschel P., Trubowitz E.: {\it Inverse spectral theory}.
Academic Press, Boston, 1987.
\bibitem [Sh] {Sh} Shen Chao-Liang:
Some inverse spectral problems for vectorial Sturm-Liouville equations. Inverse Problems {\bf
17}(5), 1253--1294 (2001).
\bibitem [SP] {SP} Samsonov B. F.; Pecheritsin A. A.:
Chains of Darboux transformations for the matrix Schrodinger equation. J. Phys. A {\bf 37}(1),
239--250 (2004).
\bibitem [Yu] {Yu} Yurko V.A.: Inverse problems for matrix Sturm-Liouville operators.
Russian J. Math. Phys. {\bf13}(1) (2006).
\end{thebibliography}
\end{document}